\documentclass{amsart} 

\usepackage{graphicx}

\usepackage{amsmath,amssymb}
\usepackage{color}   
\usepackage[all]{xy}
\usepackage{pstricks}

\newtheorem{theorem}{Theorem}[section]
\newtheorem{lemma}[theorem]{Lemma}
\newtheorem{proposition}[theorem]{Proposition}
\newtheorem{corollary}[theorem]{Corollary}

\theoremstyle{definition}
\newtheorem{definition}[theorem]{Definition}
\newtheorem{example}[theorem]{Example}

\theoremstyle{remark}
\newtheorem{remark}[theorem]{Remark}

\newcommand {\R} {\ensuremath {\mathbb{R}} }

\newcommand {\xto}[1] {\ensuremath {\xrightarrow{#1}} }
\newcommand {\xfrom}[1] {\ensuremath {\xleftarrow{#1}} }

\newcommand {\Set} {\ensuremath {\mathbf{Set}} }
\newcommand {\sSet} {\ensuremath {\mathbf{sSet}} }
\newcommand {\LPS} {\ensuremath {\mathbf{LoPospc}} }
\newcommand {\ALPS} {\ensuremath {\mathbf{A \downarrow LoPospc}} }
\newcommand {\C} {\ensuremath {\mathbf{C}} }
\newcommand {\UC} {\ensuremath {\mathbf{UC}} }
\newcommand {\UCM} {\ensuremath {\mathbf{UC/M}} }

\newcommand {\sSetCop} {\ensuremath {\sSet^{\C^{\op}}} }
\newcommand {\sPre} {\ensuremath {\mathbf{sPre}} }
\newcommand {\sPreC} {\ensuremath {\mathbf{sPre(\C)}} }
\newcommand {\sPreCOverM} {\ensuremath {\mathbf{sPre(\C) / M}} }
\newcommand {\simplCat} {\ensuremath {\mathbf{\Delta}} }

\newcommand {\dI} {\protect \overrightarrow{I}}
\newcommand {\dII} {\dI \times \dI}
\newcommand {\dO} {\protect \overrightarrow{O}}
\newcommand {\dX} {\protect \overrightarrow{X}}
\newcommand {\dS} {\protect \overrightarrow{S}}
\newcommand {\Posp} {\ensuremath {\mathbf{Pospace}} }
\newcommand {\APosp} {\ensuremath {\mathbf{A \downarrow Pospace}} }
\newcommand {\SPosp} {\ensuremath {\mathbf{S^0 \downarrow Pospace}} }
\newcommand {\from} {\leftarrow}
\newcommand {\fc} {\protect \overrightarrow{\pi}_1}

\DeclareMathOperator{\Ob}{Ob}

\DeclareMathOperator{\op}{op}

\DeclareMathOperator{\Id}{Id}

\begin{document}

\title{Context for models of concurrency} 

\author{Peter Bubenik}
\address{Institut de G\'{e}om\'{e}trie, Alg\`{e}bre et Topologie\\ 
    Ecole Polytechnique F\'{e}d\'{e}rale de Lausanne\\
    Lausanne, Switzerland} 
\email{peter.bubenik@epfl.ch}
\thanks{This work was funded by the Swiss National Science
  Foundation grant 200020-105383.}

\begin{abstract} 
Many categories have been used to model concurrency.
Using any of these, the challenge is to reduce a given model to a
smaller representation which nevertheless preserves the relevant
computer-scientific information. 
That is, one wants to replace a given model with a simpler model with
the same directed homotopy-type.
Unfortunately, the obvious definition of directed homotopy equivalence
is too coarse. 
This paper introduces the notion of \emph{context} to refine this
definition.
\end{abstract}


\date{July 15, 2005}

\maketitle

\psset{arrowsize=5pt}

\section{Introduction} \label{intro}

Various topological models are being used for studying
concurrency. 
Among them are precubical complexes~\cite{goubaultThesis},
d-spaces~\cite{grandisDHTi,grandisDHTii}, local 
po-spaces~\cite{fgrAlgebraicTaCpreprint,goubaultSomeGPiCT}, and
FLOW~\cite{gaucherModelCategory}. 
For a given concurrent system, each of these categories provides a
model which captures the relevant computer-scientific properties of the
system. 

These categories are large in two senses. 
They are large `locally' in that a given model contains many paths
which correspond to executions which are essentially equivalent. 
They are also large `globally' in that a given concurrent system has a
large number of models within the category. 
The size of these categories is a strength in terms of their
descriptive power. 
However, to aid in calculations one would like to reduce these
models to a smaller, possibly even discrete, representation.

A major goal of current research in this area is to introduce
equivalences to obtain such smaller representations, which nevertheless
still retain the relevant computer-scientific properties. 

On the local front progress has been made in reducing the path space
of a given model using directed homotopies of paths and
the fundamental category~\cite{grandisDHTi}.
One global approach is to pass to the component
category~\cite{fghrComponents,raussenStateSpaces}. 
In this paper we introduce another global approach, which is perhaps
more geometric and which is compatible with the model categorical
approach of~\cite{bubenikWorytkiewiczModelCfLPSpreprint}.

In the classical (undirected) topological case, the solution to this
`global' problem is well-understood. 
The equivalent spaces are the (weak) homotopy equivalent ones.
So for example, all of the contractible spaces (those homotopy
equivalent to a point) are equivalent. 

In the directed case there is a similar notion of directed homotopy
equivalence (abbreviated to \emph{dihomotopy equivalence}, which will
be defined in the next section). However this notion is too coarse. 

\begin{figure}[htbp]
\psset{unit=2cm}
\centering
\begin{pspicture}(1.5,1)
\psline{->}(0,0)(0.53,0.53)
\psline(0.5,0.5)(1,1)
\end{pspicture}%
\begin{pspicture}(1.5,1)
\psframe[fillstyle=solid,fillcolor=lightgray](0,0)(1,1)
\psline{->}(0.53,0)(0.5301,0)
\psline{->}(0,0.53)(0,0.5301)
\end{pspicture}%
\begin{pspicture}(1.5,1)
\psline(0.25,0)(0.75,1)
\psline(0,0.25)(1,0.75)
\psline{->}(0.4,0.3)(0.4001,0.3002)
\psline{->}(0.7,0.9)(0.7001,0.9002)
\psline{->}(0.3,0.4)(0.3002,0.4001)
\psline{->}(0.8,0.65)(0.8002,0.6501)
\end{pspicture}%
\begin{pspicture}(1,1)
\psbezier(0,0)(0,0.25)(0,0.5)(0.25,0.75)
\psbezier(0.25,0.75)(0.5,1)(0.75,1)(1,1)
\psline{->}(0.25,0.75)(0.2501,0.7501)
\pscustom{
\swapaxes
\psbezier(0,0)(0,0.25)(0,0.5)(0.25,0.75)
\psbezier(0.25,0.75)(0.5,1)(0.75,1)(1,1)
}
\psline{->}(0.75,0.25)(0.7501,0.2501)
\end{pspicture}%
\caption{$\dI$, $\dII$, $\dX$, and $\dO$} \label{figsBasic}
\end{figure}

\begin{example}
Let $\dI$ be the unit interval $[0,1]$ with a direction given by the
usual ordering of the real numbers. 
Let $\dII$ be $[0,1] \times [0,1]$ with the ordering $(x,y) \leq
(x',y')$ if and only if $x \leq x'$ and $y \leq y'$.
Let $\dX$ be the space in Figure~\ref{figsBasic} given by attaching two
copies of $\dI$ at their midpoints.
Then as will be shown explicitly in Example~\ref{egDicontractible},
$\dI$, $\dI \times \dI$ and $\dX$ are all dihomotopy equivalent to a point. 
However $\dI$ models an execution with one initial state and one final
state while $\dX$ models an execution with two initial states and two
final states.
\end{example}

Clearly a stronger notion of equivalence is needed. 
Since $\dI$ and $\dI \times \dI$ both have one initial state and one
final state and all execution paths seem to be essentially equivalent
it seems natural that we should look for a definition of equivalence under
which these are equivalent. 
However even this `equivalence' has a pitfall. 

For a notion of equivalence to be practical it should continue to hold
under certain `pastings'.
Our philosophy is the following. If we make the same addition to
equivalent models we should still have equivalent models. 

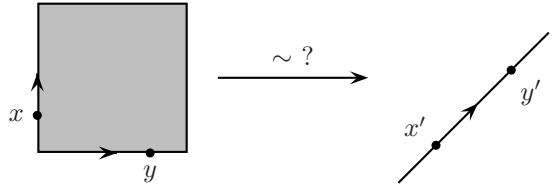
\begin{figure}[htbp]
\psset{unit=2cm}
\centering
\begin{pspicture}(0,-0.2)(1.2,1)
\psframe[fillstyle=solid,fillcolor=lightgray](0,0)(1,1)
\psline{->}(0.53,0)(0.5301,0)
\psline{->}(0,0.53)(0,0.5301)
\qdisk(0,0.25){0.03}
\qdisk(0.75,0){0.03}
\uput[l](0,0.25){$x$}
\uput[d](0.75,0){$y$}
\end{pspicture}%
\begin{pspicture}(0,-0.2)(1.2,1)
\psline{->}(0,0.5)(1,0.5)
\uput[u](0.5,0.5){$\sim$ ?}
\end{pspicture}%
\begin{pspicture}(1,1)
\psline{->}(0,0)(0.53,0.53)
\psline(0.5,0.5)(1,1)
\qdisk(0.25,0.25){0.03}
\qdisk(0.75,0.75){0.03}
\uput[ul](0.25,0.25){$x'$}
\uput[dr](0.75,0.75){$y'$}
\end{pspicture}%
\caption{A hypothetical equivalence} \label{figHypEq}
\end{figure}

\begin{example} \label{egDIIdIshouldBeNo}
Assume we have an equivalence $\dII \to \dI$ as in Figure~\ref{figHypEq}.
Consider the following pasting on $\dII$.
Let $\dO$ be the space in Figure~\ref{figsBasic} constructed by attaching
two copies of $\dI$ at their initial points and at the final
points.\footnote{This is M.Grandis' \emph{ordered circle} $\uparrow \! \!
  O^1$~\cite[Section 1.2]{grandisDHTi}.} 
Let $\dO_1$ and $\dO_2$ be two copies of $\dO$. 
For $i=1,2$ let $a_i,b_i \in \dO_i$ denote the initial and final
points of $\dO_i$.
Now choose two points $x,y \in \dII$ such that neither $x \leq y$ nor
$y \leq x$.
Let $x',y' \in \dI$ be the images of $x$ and $y$ under the assumed
equivalence (Figure~\ref{figHypEq}).
Then either $x' \leq y'$ or $y' \leq x'$, since $\dI$ is totally ordered. 

\begin{figure}
  \psset{unit=2cm}
  \centering
\begin{pspicture}(-0.75,-0.75)(1.2,1)
  \psframe[fillstyle=solid,fillcolor=lightgray](0,0)(1,1)
  \psline{->}(0.53,0)(0.5301,0) \psline{->}(0,0.53)(0,0.5301)
  \qdisk(0,0.25){0.03}
  \qdisk(0.75,0){0.03}
  \begin{pspicture}(0.5,-0.25)(0.5,-0.25)
    \psset{unit=1cm}
    \psbezier(0,0)(0.125,0.125)(0.25,0.25)(0.5,0.25)
    \psbezier(0.5,0.25)(0.75,0.25)(0.875,0.125)(1,0)
    \psline{->}(0.6,0.25)(0.6001,0.25)
    \psbezier(0,0)(0.125,-0.125)(0.25,-0.25)(0.5,-0.25)
    \psbezier(0.5,-0.25)(0.75,-0.25)(0.875,-0.125)(1,0)
    \psline{->}(0.6,-0.25)(0.6001,-0.25)
    \qdisk(0,0){0.03}
    \uput[l](0,0){$a_1$}
  \end{pspicture}%
  \begin{pspicture}(-0.75,0)(-0.75,0)
    \psset{unit=1cm}
      \psbezier(0,0)(0.125,-0.125)(0.25,-0.25)(0.25,-0.5)
      \psbezier(0.25,-0.5)(0.25,-0.75)(0.125,-0.875)(0,-1)
      \psline{->}(0.25,-0.6)(0.25,-0.6001)
      \psbezier(0,0)(-0.125,-0.125)(-0.25,-0.25)(-0.25,-0.5)
      \psbezier(-0.25,-0.5)(-0.25,-0.75)(-0.125,-0.875)(0,-1)
      \psline{->}(-0.25,-0.6)(-0.25,-0.6001)
      \qdisk(0,-1){0.03}
      \uput[d](0,-1){$b_2$}
  \end{pspicture}%
\end{pspicture}%
\begin{pspicture}(0,-0.75)(1.2,1)
  \psline{->}(0,0.5)(1,0.5) 
\end{pspicture}%
\begin{pspicture}(-0.5,-0.75)(1.5,1)
  \psline{->}(0,0)(0.53,0.53) 
  \psline(0.5,0.5)(1,1)
  \qdisk(0.25,0.25){0.03}
  \qdisk(0.75,0.75) {0.03}
  \begin{pspicture}(0.25,-0.25)(0.25,-0.25)
    \psset{unit=1cm}
    \psbezier(0,0)(0.125,0.125)(0.25,0.25)(0.5,0.25)
    \psbezier(0.5,0.25)(0.75,0.25)(0.875,0.125)(1,0)
    \psline{->}(0.6,0.25)(0.6001,0.25)
    \psbezier(0,0)(0.125,-0.125)(0.25,-0.25)(0.5,-0.25)
    \psbezier(0.5,-0.25)(0.75,-0.25)(0.875,-0.125)(1,0)
    \psline{->}(0.6,-0.25)(0.6001,-0.25)
    \qdisk(0,0){0.03}
    \uput[l](0,0){$a_1$}
  \end{pspicture}%
  \begin{pspicture}(-0.75,-0.75)(-0.75,-0.75)
    \psset{unit=1cm}
    \psbezier(0,0)(0.125,0.125)(0.25,0.25)(0.5,0.25)
    \psbezier(0.5,0.25)(0.75,0.25)(0.875,0.125)(1,0)
    \psline{->}(0.6,0.25)(0.6001,0.25)
    \psbezier(0,0)(0.125,-0.125)(0.25,-0.25)(0.5,-0.25)
    \psbezier(0.5,-0.25)(0.75,-0.25)(0.875,-0.125)(1,0)
    \psline{->}(0.6,-0.25)(0.6001,-0.25)
    \qdisk(1,0){0.03}
    \uput[r](1,0){$b_2$}
  \end{pspicture}%
\end{pspicture}%
\caption{A map $B \to C$ which should not be an equivalence} \label{figNotEq}
\end{figure}
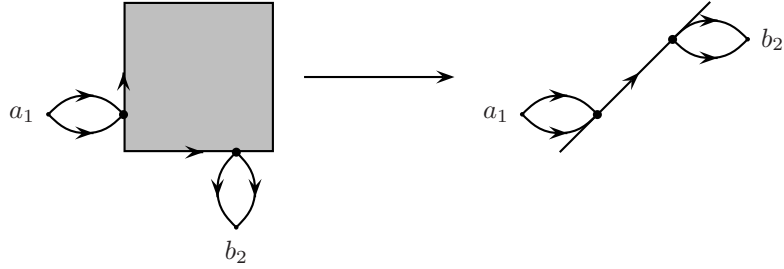

If $x' \leq y'$ then identify $b_1$ and $x$ and identify $a_2$ and $y$.
Call this space $B$ and denote $C$ the space obtained by collapsing $\dII
\subset B$ to $\dI$ using the given equivalence (Figure~\ref{figNotEq}).
Then there is an execution path from $a_1$ to $b_2$ in $C$ but not in $B$.
So the concurrent systems modeled by $B$ and $C$ are not equivalent. 
A similar construction is possible if $y' \leq x'$.
Thus from this point of view $\dII$ and $\dI$ should not be equivalent.
\end{example}
  
This gives a good indication of the current state of
affairs for determining a global notion equivalence.
We don't even know whether or not $\dII$ and $\dI$ should be equivalent. 

In this paper we introduce the idea of \emph{context}.
Whether or not $\dI$ and $\dII$ are equivalent depends on the context.
If we permit pastings as in Example~\ref{egDIIdIshouldBeNo}, then they are not
equivalent. 
However if we only permit pastings to the initial and final points of
$\dI$ and $\dII$ then they are equivalent. 
From the computer-scientific point of view this can be interpreted as
follows.
We cannot expect equivalent concurrent systems to still be
equivalent after arbitrary (but equal) changes.
However, if equal additions are made in a suitably modular way, then
the resulting systems should still be equivalent. 

It should be noted that in the examples in this paper the context is
chosen `by hand'. 
The problem of choosing the context is related to the components of
the fundamental category~\cite{fghrComponents} and to the universal
dicovering space~\cite{fajstrupDicoveringSpaces}.
A procedure for choosing the context is a subject for future research.

\noindent
\textbf{Acknowledgments.} The author would like to thank the referee
for simplifying the proof of Proposition~\ref{propFCisom}.

\section{Context for directed homotopy equivalences}

In this section we make precise the intuitive ideas presented in the
introduction.

\begin{definition}
\begin{itemize}
\item A \emph{partial order} on a topological space $U$ is a
  reflexive, transitive, anti-symmetric relation $\leq$.  If $U$ has a
  partial order $\leq$ which is a closed subset of $U \times U$ under
  the product topology, then call $U$ a \emph{po-space} or \emph{pospace}.
\item A \emph{dimap} $f: (U_1, \leq_1) \to (U_2, \leq_2)$ is a
  continuous map $f:U_1 \to U_2$ such that $x \leq_1 y$ implies that
  $f(x) \leq_2 f(y)$.
\item A product of pospaces $(U_1,\leq_1)$ and $(U_2,\leq_2)$ is a
  pospace whose underlying topological space is $U_1 \times U_2$ and
  whose order relation is given by $(x,y) \leq (x',y')$ if and only
  if $x \leq_1 x'$ and $y \leq_2 y'$.
\item A subspace $A$ of a pospace $U$ inherits a pospace structure
  under the definition $x \leq_A y$ if and only if $x \leq_U y$. This
  is called a \emph{sub-pospace}.
\end{itemize}
\end{definition}

\begin{definition} \label{defnPosp}
  Let \Posp be the the category whose objects are pospaces and whose
  morphisms are dimaps.
\end{definition}

For the sake of simplicity we will work with pospaces but one should
be able to easily extend or adapt the constructions presented here for
other models of concurrency. 

Let $\dI = ([0,1], \leq)$ where $\leq$ is the usual ordering of
$\R$. This is a pospace.
A \emph{dipath} in a pospace $B$ is a dimap $\dI \to B$.

\begin{definition} \label{defnDihomotopy}
\begin{itemize}
\item
Given dimaps $f,g: B \to C \in \Posp$, $\phi:B \times \dI \to C \in
\Posp$ is a 
\emph{dihomotopy}\footnote{This is the notion of dihomotopy 
  in~\cite{grandisDHTi} which is stronger than the notion of
  dihomotopy in~\cite{fgrAlgebraicTaCpreprint} (which uses $I=[0,1]$
  with the trivial ordering $x\leq_I y \Leftrightarrow x=y$, instead of
  $\dI$).} 
from $f$ to $g$ if $\phi|_{B \times \{0\}} = f$ and $\phi|_{B \times \{1\}} = g$.
In this case write $\phi: f \to g$. 
\item
Write $f \simeq g$ if there is a chain of dihomotopies $f \to f_1
\from f_2 \to \ldots \from f_n \to g$.
This is an equivalence relation.
\item
A dimap $f: B \to C$ is a \emph{dihomotopy equivalence} if there is
a dimap $g: C \to B$ such that $g\circ f \simeq \Id_{B}$ and $f \circ
g \simeq \Id_C$.
In this case write $B \simeq C$.
\end{itemize}
\end{definition}

Our explicit dihomotopies will often be of the following form.
\begin{definition} \label{defnLinearInterpolation}
Assume that $C$ is a pospace whose underlying topological spaces is a
subspace of $\R^n$ for some $n$.
Assume $f, g: B \to C$ are two dimaps. 
Let the \emph{linear interpolation}
between $f$ and $g$ be the map $H: B \times \dI \to \R^n$ given by
$H(b,t) = (1-t)f(b) + tg(b)$.
\end{definition}

\begin{remark}
Note that there is no guarantee that the image of such a map is in $C$. 
However one can check that it is for the cases we will consider.
\end{remark}

\begin{lemma} \label{lemmaLinearInterpolation}
Assume that $C$ is a pospace whose underlying topological spaces is a
subspace of $\R^n$ for some $n$.
If $f,g:B \to C$ are dimaps such that for all $b \in B$, $f(b) \leq
g(b)$ then if the image of the linear interpolation $H$ between $f$ and
$g$ is in $C$ then $H$ is a dihomotopy from $f$ to $g$.
\end{lemma}

\begin{proof}
That $H$ restricts to $f$ and $g$ follows from the definition of
linear interpolation.
It remains to check that $H$ is a dimap.

Let $b \leq_B b'$ and $t \leq_{\dI} t'$.
Then
\[
\begin{split}
H(b,t) &= f(b) + t(g(b) - f(b))\\
&\leq f(b) + t'(g(b)-f(b))\\
&= (1-t')f(b) + t'g(b)\\
&\leq (1-t')f(b') + t'g(b')\\
&= H(b',t')
\end{split}
\]
\end{proof}

\begin{example} \label{egDicontractible}
If $\dI_1$ and $\dI_2$ are two copies of $\dI$, then let $\dX = (\dI_1
\amalg \dI_2) / \sim$ where $(\frac{1}{2})_1 \sim (\frac{1}{2})_2$
(see Figure~\ref{figsBasic}).
We will show that under Definition~\ref{defnDihomotopy}, $\dI$,
$\dII$, and $\dX$ are dihomotopy equivalent to a point.
Let $f:\dI \to *$, $g:* \to \dI$ be the constant map and the inclusion
of the point to $1 \in \dI$.
Then $f \circ g = \Id_{*}$ and it remains to show that $\Id_{\dI}
\simeq g \circ f$.
Let $H:\dII \to \dI$ be the linear interpolation between $\Id_{\dI}$
and $g \circ f$.
That is,
\begin{eqnarray*}
H(x,t) & = & (1-t)x + t \\
& = & x + t(1-x)
\end{eqnarray*}
Then $H$ is a dimap and is the desired homotopy $\Id_{\dI} \to g \circ
f$.

In exactly the same way one can show that the constant map $f:\dII \to
*$ is a dihomotopy equivalence with $g:* \to \dII$ given by $g(*) =
(1,1)$.

To show that the constant map $f: \dX \to *$ is a dihomotopy
equivalence with $g(*) = (\frac{1}{2})_1 = (\frac{1}{2})_2$ is
slightly more complicated.
Again $f \circ g = \Id_{*}$.
To show $\Id_{\dX} \simeq g \circ f$ we will construct a chain of
dihomotopies $\Id_{\dX} \xto{H_1} h \xfrom{H_2} g \circ f$.
Let $h$ be the map that collapses the lower two line segments of $\dX$.
That is, let $h:\dX \to \dX$ be given by
\[
x \mapsto \left\{
\begin{aligned}
\frac{1}{2} & & \text{if } x < \frac{1}{2} \\
x && \text{otherwise}
\end{aligned}
\right.
\]
Let $H_1$ be the linear interpolation between $\Id_{\dX}$ and $h$ and
let $H_2$ be the linear interpolation between $g \circ f$ and $h$.
Then $H_1$ and $H_2$ are dimaps and are the desired dihomotopies.
\end{example}

We will show that in the right \emph{context} it is no longer true
that $\dI$, $\dII$, and $\dX$ are dihomotopy equivalent to a point.

\begin{definition} \label{defnContext}
Let the \emph{context} be an object $A \in \Posp$.
Instead of working in the category $\Posp$ we will work in the
category $\APosp$ of pospaces under $A$.
The objects of $\APosp$ are dimaps $A \xto{\iota_B} B$ where $B \in
\Ob \Posp$.
The morphisms in $\APosp$ are dimaps
\[
\xymatrix{
  & A \ar[dl]_{\iota_B} \ar[dr]^{\iota_C} & \\
  B \ar[rr]^{f} & & C }
\]
such that $f \circ \iota_B = \iota_C$.
\end{definition}

\begin{example}
For example if $A = S^0 = \{a,b\}$ then $B \in \Ob \APosp$ is a
pospace with two marked points.
An important example is $\dI$ with $\iota_{\dI}(a) = 0$ and
$\iota_{\dI}(b) = 1$.
\end{example}

\begin{definition} \label{dihomotopyRelA} 
\begin{itemize}
\item Given dimaps $f,g: B \to C \in \APosp$, $\phi$ is a
  \emph{dihomotopy} from $f$ to $g$ if $\phi: B \times \dI \to C \in
  \Posp$, $\phi|_{B \times \{0\}} = f$, $\phi|_{B \times \{1\}} = g$,
  and for all $a \in A$, $\phi(\iota_B(a),t) = \iota_C(a)$.  In this
  case write $\phi: f \to g$.
\item Write $f \simeq g$ if there is a chain of dihomotopies $f \to
  f_1 \from f_2 \to \ldots \from f_n \to g$.  This is an equivalence
  relation.
\item A dimap $f: B \to C$ is a \emph{dihomotopy equivalence} if there
  is a dimap $g: C \to B$ such that $g\circ f \simeq \Id_{B}$ and $f
  \circ g \simeq \Id_C$.  In this case write $B \simeq C$.
\end{itemize}
\end{definition}

We can think of this as \emph{dihomotopy rel $A$}.
In case the context $A$ is one point or two points we get pointed and
bipointed dihomotopies.
However we will see that this notion is useful for more general
contexts.

\begin{example}
Let us return to the example above.
In the context of its end points $\dI$ is no longer dihomotopic to a
point.
There is a dimap
\[
\xymatrix{
  & S^0 \ar[dl]_{\iota_{\dI}} \ar[dr]^{\iota_{*}} & \\
  {\dI} \ar[rr]^{f} & & {*} }
\]
making the diagram commute, but there is no map $g:* \to \dI$
making the diagram commute.
\end{example}

\begin{example} \label{egDIIdIyes}
In the context of $S^0 = \{a,b\}$ let $\iota_{\dI}(a) = 0$,
$\iota_{\dI}(b) = 1$,  $\iota_{\dII}(a) = (0,0)$, and
$\iota_{\dII}(b) = (1,1)$.
We claim that in this context $\dI$ and $\dII$ are dihomotopy
equivalent.
Let $f: \dII \to \dI$ and $g: \dI \to \dII$ be given by $f(x,y) =
\max(x,y)$ and $g(x) = (x,x)$. 
Then $f$ and $g$ are both dimaps, $f \circ g = \Id_{\dI}$ and $g \circ
f(x,y) = (\max(x,y), \max(x,y))$.
It remains to construct a dihomotopy rel $S^0$ from $\Id_{\dII}$ to $g
\circ f$.

Let $\phi$ be the linear interpolation (see
Definition~\ref{defnLinearInterpolation}) of $\Id_{\dII}$ and $g \circ
f$.
That is,
\begin{eqnarray*}
\phi(x,y,t) & = & (1-t)(x,y) + t(\max(x,y),\max(x,y)) \\
& = & (x + t(\max(x,y)-x), y+t(\max(x,y)-y)).
\end{eqnarray*}
Then $\phi$ is the desired dihomotopy rel $\{a,b\}$.

Hence $\dII$ and $\dI$ are dihomotopy equivalent in the given context.
\end{example}

\section{Context and the fundamental category}

We will now introduce some definitions and prove some lemmas that will
allow us to relate dihomotopy rel $A$ to the \emph{fundamental
  category}.
Furthermore it will enable us to quickly see that certain spaces are
not dihomotopy equivalent in a given context.

\begin{definition}
Let $B \in \Posp$ and let $x,y \in B$.
\begin{itemize}
\item A \emph{dipath} is a dimap $\gamma: \dI \to B$.
\item Let $\gamma_1, \gamma_2: \dI \to B$ be dipaths such that
  $\gamma_1(0) = \gamma_2(0) = x$ and   $\gamma_1(1) = \gamma_2(2) =
  y$. Then $\gamma_1$ and $\gamma_2$ are dihomotopic if they
  are dihomotopy equivalent with respect to their endpoints. That is,
  $\gamma_1 \simeq \gamma_2$ in $\SPosp$ where $\iota_{\dI}(a) = 0$,
  $\iota_{\dI}(b) = 1$,  $\iota_{B}(a) = x$, and  $\iota_{B}(b) = y$.
  In this case write $\gamma_1 \simeq \gamma_2$.
\item Let $\fc(B)(x,y)$ be the set of dihomotopy classes of
  dipaths from $x$ to $y$.
  The \emph{fundamental category} of $B$ is the category $\fc(B)$
  whose objects are the points of $B$ and whose morphisms between $x$ and $y$
  are the elements of $\fc(B)(x,y)$.\footnote{This differs from the
  definition of fundamental category in~\cite{fghrComponents}
  where the dihomotopy classes of dimaps use $I$ and not $\dI$.}
\end{itemize}
\end{definition}

\begin{lemma} \label{lemmaDihomotopyPost}
Given dihomotopic dipaths $\gamma \simeq \gamma': \dI \to B$ and a
dimap $f:B \to C$, then $f \circ \gamma \simeq f \circ \gamma'$ are
dihomotopic dipaths.
\end{lemma} 

\begin{proof}
Since $\gamma \simeq \gamma'$ there is a chain of dihomotopies $\gamma
\xto{H_1} \gamma_1 \xfrom{H_2} \gamma_2 \xto{H_3} \ldots \xfrom{H_n}
\gamma_n \xto{H_{n+1}} \gamma'$.  
Then $f \circ \gamma \xto{f \circ H_1} f \circ \gamma_1 \xfrom{f \circ H_2} f \circ \gamma_2
\xto{f \circ H_3} \ldots \xfrom{f \circ H_n} f \circ \gamma_n \xto{f
  \circ H_{n+1}} f \circ \gamma'$ is a chain of dihomotopies from $f
\circ \gamma$ to $f \circ \gamma'$. 
\end{proof}

\begin{corollary}
For a dimap $f: B \to C$ and $x,y \in B$ there is an induced map
$\fc(f): \fc(B)(x,y) \to \fc(C)(f(x),f(y))$ mapping $[\gamma] \mapsto
[f \circ \gamma]$.
That is, a dimap $f: B \to C$  induces a functor $\fc(f): \fc(B) \to \fc(C)$.
\end{corollary}

\begin{lemma} \label{lemmaDihomotopyPre}
Given dihomotopic dimaps $f \simeq g: B \to C \in \APosp$
and a dipath $\gamma: \dI \to B$ such that $\gamma(0) = \iota_B(a)$
and $\gamma(1) = \iota_B(b)$ where $a,b \in A$ then $f \circ \gamma
\simeq g \circ \gamma$ are dihomotopic dipaths.
\end{lemma}

\begin{proof}
Since $f \simeq g$ there is a chain of dihomotopies $f \xto{H_1} f_1
\xfrom{H_2} f_2 \xto{H_3} \ldots \xfrom{H_n} f_n \xto{H_{n+1}} g$. 
For $1 \leq i \leq n+1$, let $H'_i = H_i \circ (\gamma \times \dI)$.
\[
\xymatrix{
\dI \times \dI \ar[d]_{\gamma \times \dI} \ar@{-->}[dr]^{H'_i} \\
B\times \dI \ar[r]^{H_i} & C
}
\]
Then $f \circ \gamma \xto{H'_1} f_1 \circ \gamma \xfrom{H'_2} f_2 \circ \gamma
\xto{H'_3} \ldots \xfrom{H'_n} f_n \circ \gamma \xto{H'_{n+1}} g \circ
  \gamma$ is a chain of dihomotopies from the dipath $f \circ \gamma$
  to the dipath $g \circ \gamma$.
\end{proof}

\begin{proposition} \label{propFCisom}
If $f:B \to C \in \APosp$ is a dihomotopy equivalence then for all
$a,b \in A$ the induced set map $\fc(f)(a,b): \fc(B)(\iota_B(a),
\iota_B(b)) \to \fc(C)(\iota_C(a), \iota_C(b))$ is a bijection.
\end{proposition}

\begin{proof}
By definition there is a dimap $g:C \to B$ such that $g \circ f \simeq
\Id_B$ and $f \circ g \simeq \Id_C$.
So by Lemma~\ref{lemmaDihomotopyPre}, for any $a,b \in A$, any 
dipath $\gamma: \dI \to B$ such that $\gamma(0)=\iota_B(a)$ and
$\gamma(1)=\iota_B(b)$ and any dipath $\gamma': \dI \to C$ such that
$\gamma(0)=\iota_C(a)$ and $\gamma(1)=\iota_C(b)$, $g \circ f \circ
\gamma \simeq \gamma$ and $f \circ g \circ \gamma' \simeq \gamma'$. 
Hence $\fc(g)(a,b)$ is an inverse for $\fc(f)(a,b)$.
\end{proof}

\begin{example} \label{egDIIdIno}
Let $A = S^0 = \{a,b\}$ and choose any points $x,y \in \dII$ such that
$x \nleq y$ and $y \nleq x$.
Then the sets $\fc(\dII)(x,y)$ and $\fc(\dII)(y,x)$ are empty.
However for any dimap $f: \dII \to \dI$ (see Figure~\ref{figHypEq}),
either $f(x)\leq f(y)$ or $f(y) \leq f(x)$ since $I$ is totally ordered.
Furthermore one of $\fc(\dI)(f(x),f(y))$ and $\fc(\dI)(f(y),f(y))$ is nonempty.
So in the context of $\iota_{\dII}(a) = x$ and $\iota_{\dII}(b) = y$,
$\dII$ is not dihomotopy equivalent to $\dI$
since there can be no dihomotopy equivalence $f: \dII
\to \dI$ such that $\fc(f)(a,b)$ is an isomorphism.
\end{example}

\begin{example}
Let $\dX$ be the space defined earlier (see Figure~\ref{figsBasic}).
In the context of its four endpoints $(0)_1$, $(0)_2$, $(1)_1$, and
$(1)_2$, $\dX$ is not dihomotopy equivalent to $\dI$ (taking any four
not necessarily distinct points as the context for $\dI$).
Indeed, there are no dipaths from $(0)_1$ to $(0)_2$ and vice versa
(similarly for $(1)_1$ and $(1)_2$), whereas the same is not true for
the corresponding points in $\dI$.
\end{example}

\section{Finding simpler models using context}

In this section we look at two two-dimensional pospaces with a given
context.  
We show how each of them can be replaced with an equivalent
one-dimensional pospace by constructing explicit directed homotopy
equivalences. 

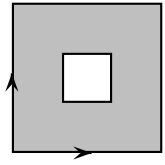
\begin{figure}[htbp]
\psset{unit=2cm}
\centering
\begin{pspicture}(1,1)
\psset{fillstyle=solid}
\psframe[fillcolor=lightgray](0,0)(1,1)
\psframe[fillcolor=white](0.3333,0.3333)(0.6666,0.6666)
\psline{->}(0.53,0)(0.5301,0)
\psline{->}(0,0.53)(0,0.5301)
\end{pspicture}%
\caption{$\dII$ with a square removed} \label{figSquareCenter}
\end{figure}

\begin{example}
In this example we show that in the context of the points $(0,0)$ and
$(1,1)$, $\dII$ with a square removed from its interior is dihomotopy
equivalent to its boundary.

Let $A = S^0 = \{a,b\}$.
Let $B$ be the sub-pospace of $\dII$ in Figure~\ref{figSquareCenter}
given by 
$ \dII \ - \ \left] \frac{1}{3},\frac{2}{3} \right[ \times
\left]  \frac{1}{3},\frac{2}{3} \right[ $.
Let $\iota_B(a) = (0,0)$ and let $\iota_B(b) = (1,1)$.
Let $C$ be the boundary of $\dII$ with  $\iota_C(a) = (0,0)$ and
$\iota_C(b) = (1,1)$. 

Intuitively we will contract $B$ to $C$ in two stages. First we will
expand the missing square $(\frac{1}{3},\frac{1}{3})\times
(\frac{1}{3},\frac{1}{3})$ to $(\frac{1}{3},1)\times(\frac{1}{3},1)$
and then to $(0,1)\times(0,1)$.
The first will be done by a map $h$ which we define below and the
composite of the two will yield the desired dihomotopy equivalence $f$.

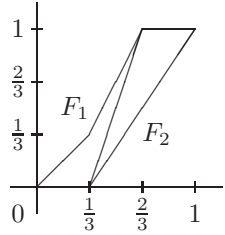
\begin{figure}[htbp]
\centering
\begin{picture}(80,80)(-10,-10)
\put(-10,0){\line(1,0){80}}
\put(0,-10){\line(0,1){80}}
\put(20,-2){\line(0,1){4}}
\put(40,-2){\line(0,1){4}}
\put(60,-2){\line(0,1){4}}
\put(-7,-10){\makebox(0,0){$0$}}
\put(-7,20){\makebox(0,0){$\frac{1}{3}$}}
\put(-7,40){\makebox(0,0){$\frac{2}{3}$}}
\put(-7,60){\makebox(0,0){$1$}}
\put(-2,20){\line(1,0){4}}
\put(-2,40){\line(1,0){4}}
\put(-2,60){\line(1,0){4}}
\put(20,-10){\makebox(0,0){$\frac{1}{3}$}}
\put(40,-10){\makebox(0,0){$\frac{2}{3}$}}
\put(60,-10){\makebox(0,0){$1$}}
\put(0,0){\line(1,1){20}}
\put(20,20){\line(1,2){20}}
\put(40,60){\line(1,0){20}}
\put(20,30){\makebox(0,0)[r]{$F_1$}}
\put(20,0){\line(2,3){40}}
\put(40,20){\makebox(0,0)[l]{$F_2$}}
\put(20,0){\line(1,3){20}}
\end{picture}%
\caption{The graphs of $F_1$, $F_2$, and $F_2 \circ F_1$.}
\label{figFgraphs}
\end{figure}

Let $F_1: [0,1] \to [0,1]$ be given by the mapping
\[ x \mapsto \left\{ 
\begin{aligned}
  x & \quad \text{if } x < \frac{1}{3} \\
  2x - \frac{1}{3} & \quad \text{if } \frac{1}{3} \leq x \leq \frac{2}{3} \\
  1 & \quad \text{if } x > \frac{2}{3}
\end{aligned}
\right.
\]
Let $F_2: [0,1] \to [0,1]$ be given by the mapping
\[ x \mapsto \left\{ 
\begin{aligned}
  0 & \quad \text{if } x < \frac{1}{3} \\
  \frac{3}{2}x - \frac{1}{2} & \quad \text{if } \frac{1}{3} \leq x
  \leq 1
\end{aligned}
\right.
\]
See Figure~\ref{figFgraphs} for graphs of $F_1$, $F_2$, and $F_2 \circ
F_1$. 

Let $f: B \to C$ and $g: C \to B$ be given by $f(x,y) = (F_2 \circ F_1(x),
F_2 \circ F_1(y))$ and  $g(x,y) = (x,y)$.
Also let $h: B \to B$ be given by $h(x,y) = (F_1(x), F_1(y))$.
Since $F_2 \circ F_1$, $\Id_{I}$, and $F_1$ are increasing maps, $f$,
$g$, and $h$ are dimaps. 

We will now give explicit dihomotopies rel $A$ showing that $g \circ f
\simeq \Id_B$ rel $A$ and $f \circ g \simeq \Id_C$ rel $A$.
Let $H_1$ and $H_2$ be linear interpolations between $\Id_B$ and $h$
and between $g \circ f$ and $h$.
That is,
\begin{gather*}
H_1(x,y,t) = (1-t)(x,y) + t(F_1(x),F_1(y)), \text{ and}\\
H_2(x,y,t) = (1-t)(F_2 \circ F_1(x),F_2 \circ F_1(y)) + t(F_1(x),F_1(y)).
\end{gather*}
Note that $F_1$ and $F_2$ fix $0$ and $1$, so $H_0$ and $H_1$ fix the
marked points $(0,0)$ and $(1,1)$.
By Lemma~\ref{lemmaLinearInterpolation}, 
$H_1: \Id_B \xto{\simeq} h$ is a dihomotopy rel $A$.
$H_2$ is a dimap since $F_2 \circ F_1(x) \leq F_1(x)$ for all $x \in
I$, so $h \leq g\circ f$.
Thus by Lemma~\ref{lemmaLinearInterpolation},
$H_2: g \circ f \xto{\simeq} h$ is a dihomotopy rel $A$.
Therefore $g \circ f \simeq \Id_B$ rel $A$ as claimed. 
Furthermore since $C$ is a sub-pospace of $B$ and $f \circ g = f = g
\circ f$, the above dihomotopies restrict to $C$ showing that $f \circ g
\simeq \Id_C$ rel $A$.
\end{example}

\begin{figure}[htbp]
\psset{unit=3cm}
\centering
\begin{pspicture}(-0.2,0)(1.2,1)
\psframe[fillstyle=solid,fillcolor=lightgray](0,0)(1,1)
\psline{->}(0.53,0)(0.5301,0)
\psline{->}(0,0.53)(0,0.5301)
\pspolygon[fillstyle=solid,fillcolor=white](0.4,0.4)(0.2,0.4)(0.2,0.6)(0.4,0.6)(0.4,0.8)(0.6,0.8)(0.6,0.6)(0.8,0.6)(0.8,0.4)(0.6,0.4)(0.6,0.2)(0.4,0.2)
\qdisk(0,0){0.02}
\qdisk(0.4,0.4){0.02}
\qdisk(0.6,0.6){0.02}
\qdisk(1,1){0.02}
\uput[l](0,0){$a$}
\uput[r](1,1){$b$}
\uput[u](0.4,0.4){$c$}
\uput[d](0.6,0.6){$d$}
\end{pspicture}%
\caption{The Swiss flag with labeled points $\{a,b,c,d\}$} \label{figSwissFlag}
\end{figure}
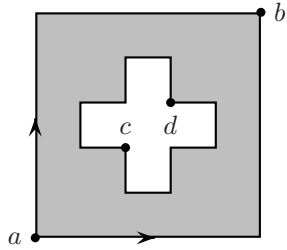
 
\begin{figure}[htbp]
\psset{unit=3cm}
\centering
\begin{pspicture}(-0.2,0)(1.2,1)
\psline(0,0)(0.4,0.4)
\psline(0.6,0.6)(1,1)
\pspolygon(0.2,0.2)(0.2,0.6)(0.4,0.8)(0.8,0.8)(0.8,0.4)(0.6,0.2)
\psline{->}(0.2,0.5)(0.2,0.5001)
\psline{->}(0.5,0.2)(0.5001,0.2)
\psline{->}(0.15,0.15)(0.1501,0.1501)
\psline{->}(0.35,0.35)(0.3501,0.3501)
\psline{->}(0.75,0.75)(0.7501,0.7501)
\psline{->}(0.92,0.92)(0.9201,0.9201)
\qdisk(0,0){0.02}
\qdisk(0.4,0.4){0.02}
\qdisk(0.6,0.6){0.02}
\qdisk(1,1){0.02}
\uput[l](0,0){$a$}
\uput[r](1,1){$b$}
\uput[u](0.4,0.4){$c$}
\uput[d](0.6,0.6){$d$}
\end{pspicture}%
\caption{A sub-pospace of the Swiss flag with the same labeled points
  $\{a,b,c,d\}$} \label{figSwissFlagSub}
\end{figure}
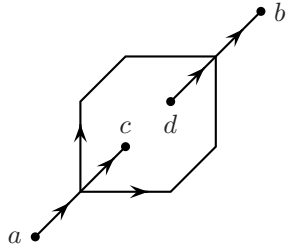

\begin{example}
The Swiss flag.

In this example we give an explicit dihomotopy between the famous
Swiss flag pospace in Figure~\ref{figSwissFlag} and the one-dimensional
sub-pospace in Figure~\ref{figSwissFlagSub} in the context of four points. 

Let $A$ be the discrete pospace $\{a,b,c,d\}$.
Let $B$ be the sub-pospace of $\dII$ given in Figure~\ref{figSwissFlag} with
the (open) cross removed and 
$\iota_B(a) = (0,0)$, $\iota_B(b) = (1,1)$, $\iota_B(c) =
(\frac{2}{5},\frac{2}{5})$, and $\iota_B(d) = (\frac{3}{5},\frac{3}{5})$. 
Let $C$ be the subspace of $B$ given in Figure~\ref{figSwissFlagSub}
with the same marked points. 

Intuitively we will contract $B$ to $C$ be applying four maps which
are described in Figure~\ref{figSwissFlagIntuitive}.

\begin{figure}
\psset{unit=3cm}
\centering
\begin{pspicture}(-0.2,0)(1.2,1)
\psframe[fillstyle=solid,fillcolor=lightgray](0,0)(1,1)
\pspolygon[fillstyle=solid,fillcolor=white](0.4,0.4)(0.2,0.4)(0.2,0.6)(0.4,0.6)(0.4,0.8)(0.6,0.8)(0.6,0.6)(0.8,0.6)(0.8,0.4)(0.6,0.4)(0.6,0.2)(0.4,0.2)
\psline(0.2,0)(0.2,1)
\psline(0,0.2)(1,0.2)
\psline(0,0)(0.2,0.2)
\psline{->}(0.03,0.15)(0.13,0.15)
\psline{->}(0.03,0.3)(0.17,0.3)
\psline{->}(0.03,0.5)(0.17,0.5)
\psline{->}(0.03,0.7)(0.17,0.7)
\psline{->}(0.03,0.9)(0.17,0.9)
\psline{->}(0.15,0.03)(0.15,0.13)
\psline{->}(0.3,0.03)(0.3,0.17)
\psline{->}(0.5,0.03)(0.5,0.17)
\psline{->}(0.7,0.03)(0.7,0.17)
\psline{->}(0.9,0.03)(0.9,0.17)
\end{pspicture}%
\begin{pspicture}(-0.2,0)(1.2,1)
\psframe[fillstyle=solid,fillcolor=lightgray](0.2,0.2)(1,1)
\psline(0,0)(0.2,0.2)
\pspolygon[fillstyle=solid,fillcolor=white](0.4,0.4)(0.2,0.4)(0.2,0.6)(0.4,0.6)(0.4,0.8)(0.6,0.8)(0.6,0.6)(0.8,0.6)(0.8,0.4)(0.6,0.4)(0.6,0.2)(0.4,0.2)
\psline(0.2,0.8)(1,0.8)
\psline(0.8,0.2)(0.8,1)
\psline(0.8,0.8)(1,1)
\psline{->}(0.3,0.97)(0.3,0.83)
\psline{->}(0.5,0.97)(0.5,0.83)
\psline{->}(0.7,0.97)(0.7,0.83)
\psline{->}(0.85,0.97)(0.85,0.87)
\psline{->}(0.97,0.3)(0.83,0.3)
\psline{->}(0.97,0.5)(0.83,0.5)
\psline{->}(0.97,0.7)(0.83,0.7)
\psline{->}(0.97,0.85)(0.87,0.85)
\end{pspicture}%
\\
\vspace{5mm}
\begin{pspicture}(-0.2,0)(1.2,1)
\psframe[fillstyle=solid,fillcolor=lightgray](0.2,0.2)(0.8,0.8)
\pspolygon[fillstyle=solid,fillcolor=white](0.4,0.4)(0.2,0.4)(0.2,0.6)(0.4,0.6)(0.4,0.8)(0.6,0.8)(0.6,0.6)(0.8,0.6)(0.8,0.4)(0.6,0.4)(0.6,0.2)(0.4,0.2)
\psline(0,0)(0.2,0.2)
\psline(0.6,0.6)(1,1)
\psline(0.2,0.6)(0.4,0.8)
\psline(0.6,0.2)(0.8,0.4)
\psline{->}(0.23,0.75)(0.33,0.75)
\psline{->}(0.63,0.75)(0.73,0.75)
\psline{->}(0.63,0.35)(0.73,0.35)
\psline{->}(0.75,0.23)(0.75,0.33)
\psline{->}(0.75,0.63)(0.75,0.73)
\psline{->}(0.35,0.63)(0.35,0.73)
\psline[linewidth=2pt]{->}(0.5,0.8)(0.65,0.8)
\psline[linewidth=2pt]{->}(0.8,0.5)(0.8,0.65)
\end{pspicture}%
\begin{pspicture}(-0.2,0)(1.2,1)
\psframe[fillstyle=solid,fillcolor=lightgray](0.2,0.2)(0.4,0.4)
\psline(0,0)(0.4,0.4)
\psline(0.6,0.6)(1,1)
\pspolygon(0.2,0.2)(0.2,0.6)(0.4,0.8)(0.8,0.8)(0.8,0.4)(0.6,0.2)
\psline{->}(0.25,0.38)(0.25,0.28)
\psline{->}(0.38,0.25)(0.28,0.25)
\psline[linewidth=2pt]{->}(0.5,0.2)(0.35,0.2)
\psline[linewidth=2pt]{->}(0.2,0.5)(0.2,0.35)
\end{pspicture}%
\caption{An intuitive sketch of the dihomotopy equivalence between the
Swiss flag and its sub-pospace}  \label{figSwissFlagIntuitive}
\end{figure}
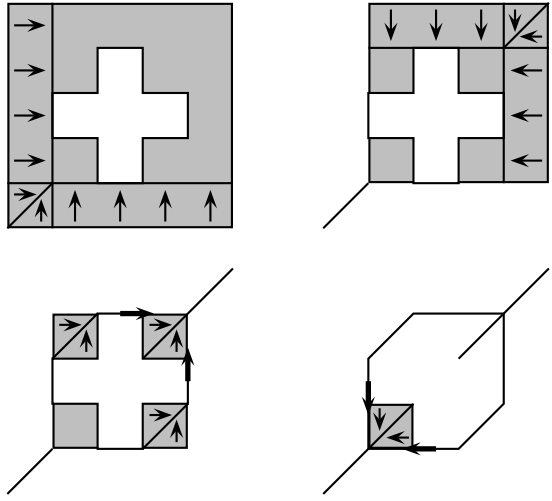

Let $g: C \to B$ be the dimap given by $g(x,y) = (x,y)$.
Let $f: B \to C$ be the dimap given by $f(x,y) = f_4 \circ f_3 \circ
f_2 \circ f_1 (x,y)$ where $f_1$, $f_2$, $f_3$, and $f_4$ are defined
in~\eqref{eqnSwissFlagfi} below.
From the sketches of $f_1,f_2,f_3$, and $f_4$ in
Figure~\ref{figSwissFlagIntuitive}, one can see that they are dimaps.
As in the previous example we will give a chain of dihomotopies rel
$A$ to show that $\Id_B \simeq g \circ f$. 
Since $C$ is a subspace of $B$ and $g \circ f = f = f \circ g$ this
will restrict to a chain of dihomotopies rel $A$ which show that
$\Id_C \simeq f \circ g$.
As a result we will have that $B \simeq C$.

\begin{gather} \label{eqnSwissFlagfi}
f_1(x,y) = \left\{ 
\begin{aligned}
  (\max(x,y), \max(x,y)) & \quad \text{if } 0 \leq x \leq \frac{1}{5},
  0
  \leq y \leq \frac{1}{5}  \\
  (\frac{1}{5}, y) & \quad \text{if } 0 \leq x \leq \frac{1}{5},
  \frac{1}{5} < y \\
  (x, \frac{1}{5}) & \quad \text{if } 0 \leq y \leq \frac{1}{5},
  \frac{1}{5} < x \\
  (x,y) & \quad \text{otherwise}
\end{aligned}
\right.
\displaybreak[0] \\
f_2(x,y) = \left\{ 
\begin{aligned} \nonumber
  (\min(x,y), \min(x,y)) & \quad \text{if } \frac{4}{5} \leq x \leq 1,
  \frac{4}{5} \leq y \leq 1 \\
  (\frac{4}{5}, y) & \quad \text{if } \frac{4}{5} \leq x \leq 1,
  y < \frac{4}{5} \\
  (x,\frac{4}{5}) & \quad \text{if } \frac{4}{5} \leq y \leq 1,
  x < \frac{4}{5} \\
  (x,y) & \quad \text{otherwise}
\end{aligned}
\right.
\displaybreak[0] \\
f_3(x,y) = \left\{ 
\begin{aligned} \nonumber
  (\max(x,y - \frac{2}{5}), \max(x+\frac{2}{5},y)) & \quad \text{if }
  \frac{1}{5} \leq x \leq \frac{2}{5}, \frac{3}{5} \leq y \leq \frac{4}{5} \\
  (\max(x,y + \frac{2}{5}), \max(x- \frac{2}{5},y)) & \quad \text{if }
  \frac{1}{5} \leq y \leq \frac{2}{5}, \frac{3}{5} \leq x \leq \frac{4}{5} \\
  (\max(x,y), \max(x,y)) & \quad \text{if }
  \frac{3}{5} \leq x \leq \frac{4}{5}, \frac{3}{5} \leq y \leq \frac{4}{5} \\
  (\frac{2}{5} + 2(x - \frac{2}{5}), y) & \quad \text{if } \frac{2}{5}
  \leq x
  \leq \frac{3}{5}, y = \frac{4}{5} \\
  (x, \frac{2}{5} + 2(y - \frac{2}{5})) & \quad \text{if }
  \frac{2}{5} \leq y
  \leq \frac{3}{5}, x = \frac{4}{5} \\
  (x,y) & \quad \text{otherwise}
\end{aligned}
\right.
\displaybreak[0] \\
f_4(x,y) = \left\{ 
\begin{aligned} \nonumber
  (\min(x,y), \min(x,y)) & \quad \text{if }
  \frac{1}{5} \leq x \leq \frac{2}{5}, \frac{1}{5} \leq y \leq \frac{2}{5} \\
  (\frac{3}{5} - 2(\frac{3}{5} - x), y) & \quad \text{if } \frac{2}{5}
  \leq x
  \leq \frac{3}{5}, y = \frac{1}{5} \\
  (x,\frac{3}{5} - 2(\frac{3}{5} - y)) & \quad \text{if } \frac{2}{5}
  \leq y
  \leq \frac{3}{5}, x = \frac{1}{5} \\
  (x,y) & \quad \text{otherwise}
\end{aligned}
\right.
\end{gather}

Let $H_1$, $H_2$, $H_3$, and $H_4$ be the linear interpolations (see
Definition~\ref{defnLinearInterpolation}) 
between $\Id_B$ and $f_1$, $f_2 \circ f_1$ and $f_1$, $f_2 \circ f_1$
and $f_3 \circ f_2 \circ f_1$, and $f$ and $f_3 \circ f_2 \circ f_1$.
Since the $f_i$ fix the labeled points, so do the $H_i$.
Furthermore, since $f_1$,$f_2$,$f_3$, and $f_4$ are dimaps, $f_1$ and $f_3$ are
increasing and $f_2$ and $f_4$ are decreasing, by
Lemma~\ref{lemmaLinearInterpolation}, the $H_i$ form  a chain of
dihomotopies 
\[
\Id_B \xto{H_1} f_1 \xfrom{H_2} f_2 \circ f_1 \xto{H_3} f_3 \circ f_2
\circ f_1 \xfrom{H_4} f = g \circ f.
\]
Therefore $\Id_B \simeq g \circ f$.
Restricting to $C$ gives a chain of dihomotopies showing $\Id_C \simeq
f = f \circ g$.
Hence $B$ is dihomotopy equivalent to $C$ rel $\{a,b,c,d\}$.
\end{example}

\section{Pushouts of dihomotopy equivalences}

In this section we elaborate on the statement made in the introduction
that dihomotopy equivalences should be preserved by `pastings'.
In fact we discuss the construction of a homotopy theory for
concurrency.
In order that we do not lose focus from the main ideas of this paper,
we will defer the details of the definitions and constructions of this
section to the appendix. 

An excellent framework for a homotopy theory on a category is given by
a \emph{model structure} on the category~\cite{hoveyBook}.
A category with a model structure and all small limits and colimits is
called a \emph{model category}.
A model structure has three special classes of morphisms:
\emph{fibrations}, \emph{cofibrations}, and \emph{weak equivalences}
which satisfy certain axioms (see Appendix~\ref{apxModelCat} for the full
definition).

The category \Posp has all small limits and colimits.
However it is too restrictive to model many concurrent systems (for
example pospaces cannot contain loops).
Though all of our examples are in \Posp a better framework for
concurrency is the category \LPS of \emph{local pospaces}.
A \emph{local pospace} is a topological space such that each point has
a neighborhood which is a pospace and that these local orders are
compatible (for a precise definition see Appendix~\ref{apxLPS}).

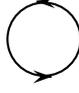
\begin{figure}[htbp]
\centering
\begin{pspicture}(1,1)
\pscircle(0.5,0.5){0.5}
\psline{->}(0.6,0)(0.6001,0)
\psline{->}(0.4001,1)(0.4,1)
\end{pspicture}
\caption{The local pospace $\dS^1$} \label{figdSone}
\end{figure}

\begin{example} \label{egdSone}
An example of a local pospace is the directed circle $\dS^1$ in
Figure~\ref{figdSone} obtained by identifying the endpoints of $\dI$.
While $\dS^1$ does not have a transitive, anti-symmetric order,
locally it has the structure of the pospace $\dI$.
\end{example}

Unfortunately, unlike $\Posp$, \LPS does not contain all small
colimits.
However there is a formal method of enlarging a category to one with
all small limits and colimits.%
\footnote{Again more details are provided in the appendix (one passes to the
  category of simplicial presheaves~\cite{duggerUHT}).}
Furthermore this larger category has a canonical model
structure!~\cite{duggerUHT}
For details on how this theory can be applied to \LPS
see the appendix and~\cite{bubenikWorytkiewiczModelCfLPSpreprint}.
In the appendix we give a more precise version of the following
theorem (Theorem~\ref{thmUCprecise}) which is proved
in~\cite{bubenikWorytkiewiczModelCfLPSpreprint}. 

\begin{theorem} \label{thmUC}
Let $\C = \LPS$. Then $\C$ is a subcategory of a model category $\UC$.
The morphisms in \C that are cofibrations are the monomorphisms and the
morphisms in \C that are weak equivalences are the isomorphisms.
\end{theorem}

From the point of view of just $\C$, this model structure is almost trivial.
However one can \emph{localize} \UC with respect to a set $M$ of
morphisms in \C to obtain a new category $\UCM$.
\UCM has the same objects and cofibrations as \UC but the morphisms in
$M$ are now weak equivalences~\cite{duggerUHT}.
The problem is to choose a good set of morphisms $M$.
For example, we can take $M$ to be the set of dihomotopy equivalences
in $\C$.

One of the key properties of \UC and \UCM is that they are \emph{left
  proper}.
That is, the pushout of a weak equivalence over a cofibration is a
  weak equivalence.
\[
\xymatrix{
  G \ar[r]_f^{\sim} \ar@{>->}[]+<0ex,-2.3ex>;[d]_(0.25)j & C
  \ar@{>->}[]+<0ex,-2.3ex>;[d] \\ 
  D \ar[r]_g^{\sim} & E }
\]
In particular in \UCM if $f \in M$ then $g$ is a weak equivalence.

\begin{example}
Recall the dihomotopy equivalence $f: \dII \to \dI$ of
Example~\ref{egDIIdIyes}.
Also recall the inclusions of $\dII$ and $\dI$ into $B$ and $C$ (see
Figure~\ref{figNotEq}) given in Example~\ref{egDIIdIshouldBeNo} where
attachments are made at the points $x,y \in \dII$ and $x',y' \in \dI$
(see Figure~\ref{figHypEq}). 
We have the following pushout diagram.
\[
\xymatrix{
  \dII \ar[r]_-f^-{\sim} \ar@{>->}[]+<0ex,-2.8ex>;[d]_(0.25)j & \dI
  \ar@{>->}[]+<0ex,-2.8ex>;[d] \\ 
  B \ar[r]_g & C }
\]
Since the inclusion $j$ is a cofibration, we get a weak equivalence
between $B$ and $C$.
However as discussed in Example~\ref{egDIIdIshouldBeNo}, $B$ should
not be equivalent to $C$. 
\end{example}

The solution to this problem is to work with \ALPS instead of \LPS
where the choice of context $A \in \Ob \LPS$ depends on the pushouts
that one would like to consider. 

In the example above the right context is clearly the points $x,y \in
\dII$ and $x',y' \in \dI$.
So $A = \{a,b\}$, $\iota_{\dII}(a) = x$, $\iota_{\dII}(b) = y$,
$\iota_{\dI}(a) = x'$, and $\iota_{\dI}(b) = y'$.
As discussed in Example~\ref{egDIIdIno} the map $f$ is not a dihomotopy
equivalence rel $A$.
So we are not forced to conclude that there is a weak equivalence
between $B$ and $C$.

In the following two examples we examine the `pastings' of two copies
of $\dII$ with a square removed. 
We show how choosing the right context allows us to find a one-dimensional
sub-pospace which is dihomotopy equivalent to the pushout.

Unlike the previous section, we will not write out the explicit dihomotopy
equivalences in these two examples.

\begin{example}
Let $A$ be the discrete space $\{a,b,c\}$.
Let $B$ be the subspace of $\dII$ in Figure~\ref{figSquareLeft} with
the square $\{(x,y) \ | \ \frac{1}{5} < x < \frac{2}{5}, \frac{2}{5} < y <
\frac{3}{5}\}$ removed.
Let $\iota_B(a) = (0,0)$, $\iota_B(b) = (\frac{1}{2},0)$, and
$\iota_B(c) = (1,0)$.

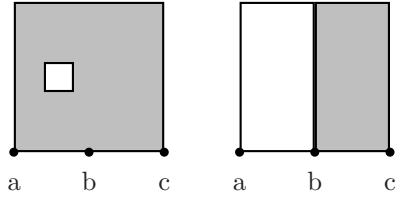
\begin{figure}
\psset{unit=2cm}
\centering
\begin{pspicture}(0,-0.25)(1.5,1)
\rput[B](0,-0.25){\psframebox[linecolor=white]{a}}
\rput[B](0.5,-0.25){\psframebox[linecolor=white]{b}}
\rput[B](1,-0.25){\psframebox[linecolor=white]{c}}
\psset{fillstyle=solid,fillcolor=lightgray}
\psframe(0,0)(1,1)
\psframe[fillcolor=white](0.2,0.4)(0.4,0.6)
\qdisk(0,0){0.03}
\qdisk(0.5,0){0.03}
\qdisk(1,0){0.03}
\end{pspicture}%
\begin{pspicture}(0,-0.25)(1,1)
\rput[B](0,-0.25){\psframebox[linecolor=white]{a}}
\rput[B](0.5,-0.25){\psframebox[linecolor=white]{b}}
\rput[B](1,-0.25){\psframebox[linecolor=white]{c}}
\psframe(0,0)(0.5,1)
\psframe[fillstyle=solid,fillcolor=lightgray](0.5,0)(1,1)
\qdisk(0,0){0.03}
\qdisk(0.5,0){0.03}
\qdisk(1,0){0.03}
\end{pspicture}%
\caption{The spaces $B$ and $B'$, which are subspaces of $\dII$ with a
  rectangle removed and labeled points $\{a,b,c\}$} \label{figSquareLeft}
\end{figure}

Let $C$ be the subspace of $\dII$ in Figure~\ref{figSquareRight} with
the square $\{(x,y) \ | \ \frac{3}{5} < x < \frac{4}{5}, \frac{2}{5} < y <
\frac{3}{5}\}$ removed.
Let $\iota_B(a) = (0,1)$, $\iota_B(b) = (\frac{1}{2},1)$, and
$\iota_B(c) = (1,1)$.

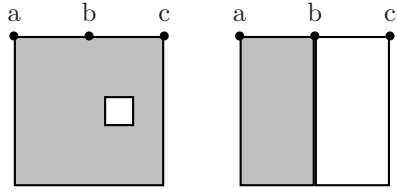
\begin{figure}
\psset{unit=2cm}
\centering
\begin{pspicture}(1.5,1.25)
\rput[B](0,1.1){\psframebox[linecolor=white]{a}}
\rput[B](0.5,1.1){\psframebox[linecolor=white]{b}}
\rput[B](1,1.1){\psframebox[linecolor=white]{c}}
\psset{fillstyle=solid,fillcolor=lightgray}
\psframe(0,0)(1,1)
\psframe[fillcolor=white](0.6,0.4)(0.8,0.6)
\qdisk(0,1){0.03}
\qdisk(0.5,1){0.03}
\qdisk(1,1){0.03}
\end{pspicture}%
\begin{pspicture}(1,1.25)
\rput[B](0,1.1){\psframebox[linecolor=white]{a}}
\rput[B](0.5,1.1){\psframebox[linecolor=white]{b}}
\rput[B](1,1.1){\psframebox[linecolor=white]{c}}
\psframe(0.5,0)(1,1)
\psframe[fillstyle=solid,fillcolor=lightgray](0,0)(0.5,1)
\qdisk(0,1){0.03}
\qdisk(0.5,1){0.03}
\qdisk(1,1){0.03}
\end{pspicture}%
\caption{The spaces $C$ and $C'$, which are subspaces of $\dII$ with a square
removed and  labeled points $\{a,b,c\}$} \label{figSquareRight}
\end{figure}

Let $B'$ be the subspace of $\dII$ in Figure~\ref{figSquareLeft} with
the rectangle 
$\left] 0, \frac{1}{2} \right[ \times \left] 0, 1 \right[$
removed and the same marked points.
Then there is a dihomotopy equivalence $f:B \xto{\simeq} B'$ rel $A$.
One can construct the required dihomotopies by stretching the region
$\frac{2}{5} \leq y \leq \frac{3}{5}$ first to $y=1$ and then to
$y=0$.
Next one stretches the region $\frac{1}{5} \leq x \leq \frac{2}{5}$
first to $x = \frac{1}{2}$ and then to $x=0$.
All this is done while leaving the three marked points fixed.

Similarly there is a dihomotopy equivalence $g:C \xto{\simeq} C'$ rel
$A$ where $C'$ is the subspace of $\dII$ in
Figure~\ref{figSquareRight} with the rectangle $\{(x,y) | \frac{1}{2} < x
< 1, 0 < y < 1\}$ removed.  

Let $D$ be the space obtained by attaching $B$ along its bottom edge
to the top edge of $C$.
Notice that $D \in \Ob \APosp$ and the inclusions $i:B \to D$ and $j:C
\to D$ are dimaps in $\APosp$.

Now take the following pushout.
\[
\xymatrix{
  B \ar[r]_f^{\sim} \ar@{>->}[]+<0ex,-2.3ex>;[d]_(0.25)i & B'
  \ar@{>->}[]+<0ex,-2.3ex>;[d] \\ 
  D \ar[r]_{f'}^{\sim} & E }
\]
Then $E$ is the pospace obtained by attaching the bottom edge of $B'$
to the top edge of $C$.
Since $C$ includes into $E$ we can take the following pushout.
\[
\xymatrix{
  C \ar[r]_g^{\sim} \ar@{>->}[]+<0ex,-2.3ex>;[d]_(0.25)i & C'
  \ar@{>->}[]+<0ex,-2.3ex>;[d] \\ 
  E \ar[r]_{f''}^{\sim} & F }
\]
Now $F$ is the pospace\footnote{Being precise, if we consider the
  pushout $F$ to be a subspace of $I \times I$ then the points
  $[0,\frac{1}{5}]\times\{\frac{1}{2}\}$ are identified as are the
  points $[\frac{2}{5},\frac{3}{5}]\times\{\frac{1}{2}\}$ and the points
  $[\frac{4}{5},1]\times\{\frac{1}{2}\}$. However this pospace is
  dihomotopy equivalent to the pospace obtained by attaching the
  bottom edge of $B'$ to the top edge of $C'$.} \label{footLeftRight}
in Figure~\ref{figSquaresLeftRight} obtained by
attaching the bottom edge of $B'$ to the top edge of $C'$.

\begin{figure}
\centering
\begin{pspicture}(2,2)
\psframe(0,0)(1,2)
\psframe[fillstyle=solid,fillcolor=lightgray](0,0)(0.5,1)
\psframe[fillstyle=solid,fillcolor=lightgray](0.5,1)(1,2)
\qdisk(0,1){0.05}
\qdisk(0.5,1){0.05}
\qdisk(1,1){0.05}
\psline{->}(0.4,0)(0.4001,0)
\psline{->}(0,0.55)(0,0.5501)
\end{pspicture}%
\psset{unit=2cm}
\begin{pspicture}(-0.2,0)(1,1)
\psframe(0,0)(1,1)
\psline(0,0)(1,1)
\qdisk(0,0){0.03}
\qdisk(0.5,0.5){0.03}
\qdisk(1,1){0.03}
\uput[l](0,0){$a$}
\uput[ul](0.5,0.5){$b$}
\uput[r](1,1){$c$}
\psline{->}(0.55,0)(0.5501,0)
\psline{->}(0,0.55)(0,0.5501)
\psline{->}(0.35,0.35)(0.3501,0.3501)
\end{pspicture}
\caption{The pospaces $F$ and $G$} \label{figSquaresLeftRight}
\end{figure}
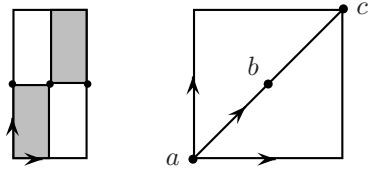

Finally $F$ is dihomotopy equivalent rel $A$ to the space $G$ in
Figure~\ref{figSquaresLeftRight}.
Consider $F$ and $G$ as sub-pospaces of $\dII$.
The dihomotopy is obtained by first collapsing the square
$[\frac{1}{2},1] \times [\frac{1}{2},1]$ using $(x,y) \mapsto
(\max(x,y), \max(x,y))$, and then collapsing the square
$[0,\frac{1}{2}] \times [0,\frac{1}{2}]$ using $(x,y) \mapsto
(\min(x,y), \min(x,y))$

Thus in the context of $A$, $D$ is equivalent to $G$.
\end{example}

\begin{example}
Let $A$, $B$, $C$, $B'$ and $C'$ be as in the previous example,
except that the marked points on $B$ and $B'$ are taken to be on the
top edge, and the marked points on $C$ and $C'$ are taken to be on the
bottom edge.
Let $D'$  be the space obtained by attaching $C$ along its bottom edge
to the top edge of $B$.

Then as in the previous example $D'$ is dihomotopy equivalent to $F'$ where
$F'$ is the pospace in Figure~\ref{figSquaresRightLeft} obtained by
attaching the bottom edge of $C'$ to the top edge of $B'$.

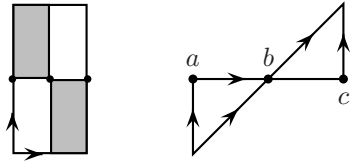
\begin{figure}
\centering
\begin{pspicture}(2,2)
\psframe(0,0)(1,2)
\psframe[fillstyle=solid,fillcolor=lightgray](0,1)(0.5,2)
\psframe[fillstyle=solid,fillcolor=lightgray](0.5,0)(1,1)
\qdisk(0,1){0.05}
\qdisk(0.5,1){0.05}
\qdisk(1,1){0.05}
\psline{->}(0.4,0)(0.4001,0)
\psline{->}(0,0.55)(0,0.5501)
\end{pspicture}%
\psset{unit=2cm}
\begin{pspicture}(-0.2,0)(1,1)
\pspolygon(0,0.5)(1,0.5)(1,1)(0,0)
\qdisk(0,0.5){0.03}
\qdisk(0.5,0.5){0.03}
\qdisk(1,0.5){0.03}
\uput[u](0,0.5){$a$}
\uput[u](0.5,0.5){$b$}
\uput[d](1,0.5){$c$}
\psline{->}(0,0.3)(0,0.3001)
\psline{->}(1,0.8)(1,0.8001)
\psline{->}(0.35,0.5)(0.3501,0.5)
\psline{->}(0.3,0.3)(0.3001,0.3001)
\psline{->}(0.8,0.8)(0.8001,0.8001)
\end{pspicture}
\caption{The pospaces $F'$ and $G'$} \label{figSquaresRightLeft}
\end{figure}

Finally $F'$ is dihomotopy equivalent rel $A$ to the space
$G'$ in Figure~\ref{figSquaresRightLeft}.
Consider $F'$ and $G'$ as sub-pospaces of $\dII$.
The dihomotopy is obtained by collapsing the regions 
$[\frac{1}{2},1] \times [0,\frac{1}{2}]$ using $(x,y) \mapsto
(x, \frac{1}{2})$, and then collapsing the square
$[0,\frac{1}{2}] \times [\frac{1}{2},1]$ using $(x,y) \mapsto
(x, \frac{1}{2})$.

Thus in the context of $A$, $D'$ is equivalent to $G'$.
\end{example}

\begin{example} \label{egNonDiscreteContext}
Finally we give an example which requires a non-discrete context.
Let $X = \dII$. We will show that if we want to use $X$ to construct a
certain space $Z$ then there is no appropriate finite context.

Let $I$ be the unit interval $[0,1]$ together with the trivial
partial order given by $x \leq_I y$ if and only if $x=y$.

\begin{figure}[htbp]
\centering
\begin{pspicture}(1.5,1)
\psframe[fillstyle=solid,fillcolor=lightgray](0,0)(1,1)
\psline{->}(0.53,0)(0.5301,0)
\psline{->}(0,0.53)(0,0.5301)
\psline(0,1)(1,0)
\end{pspicture}%
\begin{pspicture}(1,1)
\psframe[hatchangle=0,hatchwidth=0.2pt,hatchsep=2pt,fillstyle=vlines*,fillcolor=lightgray](0,0)(1,1)
\psline{->}(0,0.4)(0,0.4001)
\psline{->}(1,0.4)(1,0.4001)
\psline(0,0.5)(1,0.5)
\end{pspicture}%
\caption{$X$ and $Y$ with the images of $I$ marked} \label{figXYI}
\end{figure}
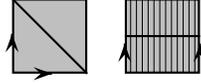

Let $\varphi: I \to X$ be the inclusion of the anti-diagonal, given by
$t \mapsto (t,1-t)$ (see Figure~\ref{figXYI}).
Let $Y = I \times \dI$ and let $\psi: I \to Y$ be the inclusion of the
central line, given by $t \mapsto (t, \frac{1}{2})$ (see
Figure~\ref{figXYI}). 
Define the pospace $Z$ obtained by gluing $X$ and $Y$ together along
the images of $I$. That is, $Z$ is the following pushout.
\[
\xymatrix{
I \ar[r]^{\varphi} \ar@{>->}[]+<0ex,-2.3ex>;[d]_{\psi} & X
\ar@{>->}[]+<0ex,-2.3ex>;[d]^{\iota_X} \\ 
Y \ar[r]_{\iota_Y} & Z
}
\]
We claim that if we want to consider this pushout then there is no
appropriate finite context.

For $\alpha \in I$ let $p_{\alpha} := \iota_X (\varphi (\alpha)) =
\iota_Y(\alpha,\frac{1}{2})$, $p^0_{\alpha} = \iota_Y (\alpha,0)$ and
$p^1_{\alpha} = \iota_Y (\alpha,1)$.
Notice that for $s \neq t \in I$ there does not exist a dipath in $Z$
from $p^0_s$ to $p^1_t$.

Now let $A$ be some context and fix $\iota_I:A \to I$ which determines
$\iota_X:A \to X$ and $\iota_Y:A\to Y$. 
Let $f: X \to X'$ be some dihomotopy equivalence rel $A$. 
Let $Z'$ and $g$ be defined by the following pushout. 
\[
\xymatrix{
X \ar[r]_{f}^{\simeq} \ar@{>->}[]+<0ex,-2.3ex>;[d]_{\iota_X} & X'
\ar@{>->}[]+<0ex,-2.3ex>;[d] \\ 
Z \ar[r]_{g}^{\simeq} & Z'
}
\]

Assume there exists $s \neq t \in I$ such that $f(\varphi(s)) = f(\varphi(t))$.
We claim that there is a dipath from $g(p^0_s)$ to $g(p^1_t)$.
In $Z$ there is a dipath from $p^0_s$ to $p_s$ and a
dipath from $p_t$ to $p^1_t$.
The concatenation of the images of these paths under $g$ gives the
desired dipath in $Z'$. 
But this contradicts Proposition~\ref{propFCisom}.

Therefore there should not have been an equivalence $f$ such that
$f(\varphi(s)) = f(\varphi(t))$ for some $s \neq t \in I$. 
We can prevent this difficulty if we use the context $A=I$ together
with $\iota_I = \Id_I$.

On the other hand with any finite context $A$, we claim that there is
always such a dihomotopy equivalence $f$.
Assume that $A$ is a finite context and fix $\iota_I:A \to I$.
Since $A$ is finite there is some interval $[a,b] \subset I \backslash
\iota_I(A)$.
Let $f:X \to X$ be the dihomotopy equivalence given by the
concatenation of the following two maps.
First collapse the region $[a,b] \times I$ to the right.
Then collapse the region $I \times [1-b,1-a]$ upwards.
Then $f$ is a dihomotopy equivalence rel $A$ but $f(\varphi(a))=f(\varphi(b))$.
\end{example}

\noindent
\textbf{Acknowledgments.} 
I would like to thank Eric Goubault, Kathryn Hess, Krzysztof
Worytkiewicz, and Emmanuel Haucourt for introducing me to the study of
concurrency using topology and category theory and for many helpful
discussions.

\appendix

\section{Model Categories} \label{apxModelCat}

In this section we define model categories, and show how a given small
category can be embedded into a \emph{universal model category}.
For more details see~\cite{duggerUHT,bubenikWorytkiewiczModelCfLPSpreprint}.

\begin{definition}
A \emph{model category} is a category $\C$ with three
distinguished classes of morphisms: weak equivalences, cofibrations,
and fibrations satisfying the following conditions:
\begin{enumerate}
\item \C contains all small limits and colimits.
\item If there exist morphisms $f$, $g$ and $g \circ f$ and two of
  them are weak equivalences then so is the third.
\item Weak equivalences, cofibrations, and fibrations are closed under
  retracts.
\item Given any commutative diagram
\[
\xymatrix{
A \ar[r] \ar[d]_i & X \ar[d]^p \\
B \ar[r] & Y
}
\]
such that $i$ is a cofibration and $p$ is a fibration, then if either
$i$ or $p$ is also a weak equivalence then there exists a map $B \to
X$ making the diagram commute. 
\item Any map may be factored as a cofibration followed by a fibration
  which is a weak equivalence, and as a cofibration which is a weak
  equivalence followed by a fibration.
\end{enumerate}
\end{definition}

Next we define the category of simplicial presheaves.

\begin{definition} \label{defnSimplicialPresheaves}
\begin{itemize}
\item The simplicial category $\simplCat$ is the category whose
  objects are $[n] = \{0,1,\ldots,n\}$ for $n\geq 0$ and whose
  morphisms are maps $f:[n] \to [k]$ such that $x \leq y$ implies that
  $f(x) \leq f(y)$.  
\item The category of simplicial sets \sSet is the category
  $\mathbf{\Set^{{\simplCat}^{\op}}}$ whose objects are contravariant
  functors from $\simplCat$ to the category of sets $\Set$ and whose
  morphisms are natural transformations.
\item Let \C be a small category. Then \sPreC is the category \sSetCop
  whose objects are the contravariant functors from \C to \sSet and
  whose morphisms are natural transformations.
\end{itemize}
\end{definition}

\begin{remark}
An important fact is that there is an embedding $\C \to \sPreC$.
\end{remark}

The category \sSet has a model structure in which the cofibrations are
the monomorphisms and the weak equivalences are the morphisms $f$ such
that $|f|$ the geometric realization of $f$ is a weak equivalence in
the category of topological spaces (that is, it induces isomorphisms
between homotopy groups). For more details see~\cite{hoveyBook}.

The category of simplicial presheaves has a canonical model structure,
called the \emph{cofibrant model structure}, where the weak
equivalences and the cofibrations are defined objectwise.
That is, a morphism $f$ in $\sPreC$ is a weak equivalence or cofibration
if and only if for each $X \in \Ob \C$ the morphism $f(X)$ is a weak
equivalence or cofibration in $\sSet$.

Now one can localize this model category~\cite{hirschhornBook} with
respect some set of morphisms $M$ to get a new model category $\sPreCOverM$.
This model category has the same objects, but in addition to the
previous weak equivalences, the morphisms in $M$ are now weak
equivalences.
For example if $\C = \LPS$ then one could localize with respect to all
dihomotopy equivalences (it makes sense to say this because of the
embedding of \C in $\sPreC$).

\section{Local po-spaces} \label{apxLPS}

In this section we give a precise definition of the category \LPS of
local pospaces and use it to give a more precise version of
Theorem~\ref{thmUC}. Local pospaces are defined
in~\cite{fgrAlgebraicTaCpreprint,bubenikWorytkiewiczModelCfLPSpreprint}. Here we
follow~\cite{bubenikWorytkiewiczModelCfLPSpreprint}. 

\begin{definition} \label{defnOrderAtlas}
  \begin{itemize}
  \item Given a topological space $M$, an \emph{order atlas} on $M$ is
    an open cover\footnote{That is, each $U_i$ is an open subset of
    $M$, and $M = \cup_{i \in I} U_i$.} $U = \{U_i\}$ indexed by a
    set $I$ such that each $U_i$ is a pospace and that the orders are
    compatible.  That is, given $x,y \in U_i \cap U_j$, $x \leq_i y$ if
    and only if $x \leq_j y$.
  \item Let $U = \{U_i\}$ and $V = \{V_j\}$ be two order atlases. Then
    $V$ is said to be a \emph{refinement} of $U$ if for any $U_i$ and
    any $x \in U_i$ there exists a $V_j$ containing $x$ which is a
    sub-pospace of $U_i$.
  \item Two order atlases are said to be \emph{equivalent} if they have a
    common refinement. One can check that this defines an equivalence
    relation.
  \item Define a \emph{local pospace} to be a topological space
    together with an equivalence class of order atlases.
  \item Define a \emph{dimap of local pospaces} $f: (M,\bar{U}) \to
    (N,\bar{V})$ to be a continuous map $f:M \to N$ such that for any
    choice of $V = \{V_j\} \in \bar{V}$ there is some choice of $U =
    \{U_i\} \in \bar{U}$ such that for all $i,j$ the partial map $f:
    U_i \to V_j$ is a dimap of pospaces.
  \end{itemize}
\end{definition}

\begin{definition} \label{defnLPS}
  Define \LPS to be the category whose objects are local pospaces whose
  underlying topological spaces are subsets of $\R^n$ for some
  $n$,\footnote{The local partial order need not be the one
    inherited from the usual partial order on $\R^n$.}  
  and whose morphisms are dimaps between local pospaces.  
\end{definition}

\begin{remark}
Notice that we have restricted the class of local pospaces in our
category. This is done precisely so that the resulting category \LPS
is a small category, 
which is used to apply the machinery of Appendix~\ref{apxModelCat}.
For the purposes of concurrency, this does not seem to be a significant
limitation.
Furthermore, it may be possible that any local pospace can be `found' in
$\sPre(\LPS)$.

Nevertheless, a consequence of this, is that the category \Posp in
Definition~\ref{defnPosp} is not a subcategory of $\LPS$.
Of course one could define a new category $\mathbf{Pospace'}$ whose
objects are those pospaces whose underlying topological spaces are
subsets of $\R^n$ for some $n$. Then $\mathbf{Pospace'}$ is a
subcategory of $\LPS$.
All of our examples are in $\mathbf{Pospace'}$.
\end{remark}

We can now give a more precise version of Theorem~\ref{thmUC}.

\begin{theorem}[\cite{bubenikWorytkiewiczModelCfLPSpreprint}] \label{thmUCprecise}
There exists a model structure on $\sPre(\LPS)$ such that the
cofibrations are the monomorphisms. 
Furthermore the morphisms in $\LPS$ which are weak equivalences in
$\sPre(\LPS)$ are just the isomorphisms.
\end{theorem}


\end{document}